\documentclass[smallextended]{svjour3}       
\smartqed  
\usepackage{graphicx}
\usepackage{mathptmx}      
\usepackage[utf8]{inputenc}
\usepackage{caption}
\usepackage{amsmath}
\usepackage{amsfonts}
\usepackage{fancybox}
\usepackage[hyphens]{url}
\usepackage[utf8]{inputenc}
\usepackage{xspace}
\usepackage{dcolumn}
\usepackage{rccol}
\usepackage{ifpdf}
\usepackage[pdftex,bookmarks=true]{hyperref}

\let\phi=\varphi

\newcommand{\systemname}[1]{\textsf{#1}\xspace}
\newcommand{\libname}[1]{\textrm{#1}\xspace}
\newcommand{\mizar}{\systemname{Mizar}}

\newcommand{\mml}{\libname{MML}}

\newcommand{\theory}[1]{\textrm{#1}}
\newcommand{\ZF}[0]{\theory{ZF}}
\newcommand{\ZFC}[0]{\theory{ZFC}}

\newcommand{\TG}[0]{\theory{TG}}

\newcommand{\N}{\mathbf{N}}
\newcommand{\equivalence}[2]{{#1} \leftrightarrow {#2}}
\newcommand{\coloneqq}[2]{{#1} \mathrel{\mathop:}= {#2}}

\newcommand{\battery}[0]{\{\forall,\wedge,\neg\}}
\newcommand{\proves}[2]{{#1} \vdash {#2}}
\usepackage{xcolor}

\newtheorem{fact}{Fact}

\begin{document}
\title{Sentence complexity of theorems in {\mizar}}
\titlerunning{Sentence complexity of theorems in {\mizar}}
\author{Jesse Alama}
\institute{Jesse Alama{\at}Center for Artificial Intelligence{\\}New University of Lisbon{\\}\email{j.alama@fct.unl.pt}}
\date{Received: date / Accepted: date}
\maketitle
\begin{abstract}{As one of the longest-running computer-assisted formal mathematics projects, large tracts of mathematical knowledge have been formalized with the help of the {\mizar} system.  Because {\mizar} is based on first-order classical logic and set theory, and because of its emphasis on pure mathematics, the {\mizar} library offers a cornucopia for the researcher interested in foundations of mathematics.  With {\mizar}, one can adopt an experimental approach and take on problems in foundations, at least those which are amenable to such experimentation.  Addressing a question posed by H.~Friedman, we use {\mizar} to take on the question of surveying the sentence complexity (measured by quantifier alternation) of mathematical theorems.  We find, as Friedman suggests, that the sentence complexity of most {\mizar} theorems is universal ($\Pi_{1}$, or $\forall$), and as one goes higher in the sentence complexity hierarchy the number of {\mizar} theorems having these complexities decreases rapidly.  The results support the intuitive idea that mathematical statements, even when carried out an abstract set-theoretical style, are usually quite low in the sentence complexity hierarchy (not more complex than $\forall\exists\forall$ or $\exists\forall$.)}

\end{abstract}
\section{Introduction}

{The formalization of mathematical reasoning has long been a venerated activity, even before the invention of modern-day computers.  As is well known, the use of computers to assist with formalized mathematical reasoning is quite old~\cite{davis2001history}.  Early on in this history, proof search (the problem of determining whether a conclusion follows from a set of premises) was separated from proof checking (the problem of determining whether a formal derivation is in fact a derivation).  In the 1960s, for instance, H.~Wang emphasized proof checking as an area where computers could potentially be of great use:}

\begin{quote}
While formalizing known or conjectured proofs and proving new theorems are intimately related, it is reasonable to suppose that the first type of problem is much easier for the machine. That is why the writer believes that machines may become of practical use more quickly for mathematical research, not by proving new theorems, but by formalizing and checking outlines of proofs. This proof formalization could be developed, say, from textbooks to detailed formulations more rigorous than \emph{Principia} [\emph{Mathematica}], from technical papers to textbooks, or from abstracts to technical papers.~\cite{wang1960toward}
\end{quote}

{Then, in the 1970s, the {\mizar} system was born.\footnote{It is not necessary to go into more detail about the early history of {\mizar} and its place in the history of automated reasoning.  The present article is part of a celebration of {\mizar}'s first $40$ years.  For an authoritative discussion of the first $30$ years, see~\cite{mizar-the-first-30-years}.  For the present author, a major feature of {\mizar} added in the last decade is the ``XMLization'' of {\mizar} begun by Urban~\cite{DBLP:conf/mkm/Urban05}, and which continues to the current day (e.g.,~\cite{DBLP:conf/aisc/2012}).  The XMLization of various aspects of {\mizar} has opened up the system to experimenters who, like the author, wish to get their hands dirty with {\mizar} but who do not have enough familiarity with the {\mizar} codebase to confidently modify the internal workings of the system.}}

{{\mizar} is a language and a suite of tools for assisting with the formalization of mathematics quite in line with what Wang envisioned.  Going even further, one can see {\mizar} as a platform for Wang's inferential analysis with an eye towards traditional foundations of mathematics (FOM).  The present paper is a contribution along these lines.}

{In more detail: {\mizar} is based on Tarski-Grothendieck set theory ($\TG$)~\cite{TARSKI.ABS,tarski1939wellordered}.  $\TG$ is stronger than the standard Zermelo-Fraenkel set theory $\ZF$ (even $\ZFC$, in which one includes the Axiom of Choice) because in $\TG$ one can prove the existence, for any set $S$, of the so-called universe for $S$, which is a set $U$ containing $S$ closed under taking subsets and power sets and which also satisfies the condition:}
\[
\text{If $X \subseteq U$ and $|X| < |U|$, then $X \in U$}.
\]
{The universe axiom of $\TG$ yields the existence of strongly inaccessible cardinal numbers, a principle that is independent of $\ZF$ (as well as $\ZFC$)~\cite{kunen2011set,jech2008axiom}.}

{Thanks to the rather long history of {\mizar} and the steady guidance of its founder, Andrzej Trybulec, a rather large array of matematical theorems can be found in the {\mizar} Mathematical Library {\mml}~\cite{Trybulec:2013:FMM}.  One could say that more or less an entire undergraduate mathematics curriculum, as well as a number of more advanced results, have been formalized in {\mizar}.}

{From the perspective of foundations of mathematics, the {\mml} can be viewed as a large collection of theorems proved in set theory.  The {\mml} can serve as the source for many potential experiments of interest for research in foundations of mathematics.  Some experiements along these lines have already been carried out~\cite{Pak:2013:MLE,GrabowskiTAC:2012,DBLP:conf/cade/Urban07,DBLP:conf/mkm/Urban03,DBLP:journals/jar/Alama13,alama-mamane-urban}.}

{Concerning the complexity of theorems of mathematics, one notices that usually mathematicians avoid many alternating quantifiers, that is, universals within existentials within universals\dots .  Some theorems, or even definitions, seem to be somewhat inherently complex in terms of quantifier alternation, such as the definition of a differentiable function or Borel set.  Commenting on the possibility of reformulating apparently difficult (many quantifier alternations) mathematical statements in more compact ways, H.~Friedman conjectured:}
\begin{quote}
Call a sentence of predicate calculus logically tight if and only if it is not logically equivalent to any sentence of lower complexity. The conjecture is that sentences occurring in actual mathematics are (generally) logically tight. This conjecture can be interpreted in two different ways. Firstly, without unraveling abbreviations. This seems like a safe bet since mathematical statements seem to always be quite short. Once they threaten to get long, one brings in abbreviations. Secondly, with unraveling abbreviations. Then the conjecture seems much more problematic and needs systematic investigation.~\cite{friedman1998complexity}
\end{quote}

{One step toward addressing this conjecture is to first of all get a sense of the actual sentence complexity of a wide range of formalized mathematical theorems.  In this note we report on a experiment conducted with the {\mml} to measure the sentence complexity (defined in Section~2) of the theorems in {\mizar}.  Section~3 reports on the data, which is discussed further in Section~4.  Section~5 proposes further problems and concludes.}

\section{Sentence complexity of {\mizar} theorems}

{To measure the complexity of mathematical theorems, we need a metric that allows us to say that one mathematical statement is more complex than another.  We are interested, in the first place, in measuring quantifier alternation.  To accomplish this, we use the familiar notion of prenex normal form, with the help of which we can define a hierarchy of formulas.}

\begin{definition}[Sentence complexity hierarchy]
{Given a first-order signature $\sigma$, define $\Sigma_{0} = \Pi_{0}$ as the class of $\sigma$-atomic formulas.  For a positive natural number $n$, define}
\[
\Sigma_{n} \coloneqq \{ \exists x_{1}, \dots, x_{k} \phi \colon \phi \in \Pi_{n - 1} \} {\quad} \text{and} {\quad} \Pi_{n} = \coloneqq \{ \forall x_{1}, \dots, x_{k} \phi \colon \phi \in \Pi_{n - 1} \}
\]
A formula $\phi$ is in prenex normal form if there exists a natural number $k$ such that $\phi \in \Pi_{k}$ or $\phi \in \Sigma_{k}$.
\end{definition}
{Thus, from an atomic formula $\phi$ one obtains a $\Pi_{1}$ formula by prefixing $\phi$ by any (positive) number of universal quantifiers, and similarly by prefixing $\phi$ by a some existential quantifiers, one obtains a $\Sigma_{1}$ formula.}

\begin{fact}
{For every formula $\phi$ there exists a formula $\phi^{*}$ in prenex normal form such that $\vdash{\phi} \leftrightarrow {\phi^{*}}$.}
\end{fact}

{Our measure of sentence complexity is entirely syntactic.  To compute the complexity of a formula, one does not need to refer to provability relative to a theory such as Peano Arithmetic or $\ZFC$, nor does one need to use structures (such as truth of a sentence in the standard structure~$\N$ of the natural numbers, as in the notion of arithmetical hierarchy~\cite{enderton2001mathematical}).  In practice, the construction of a prenex normal form, or even the construction of all possible prenex normal forms of a given formula, is clearly bounded by the complexity of $\phi$ (can count the number of subformulas of $\phi$ for which one of the prenex normal form rewrite operations applies).}

{Prenex normal form is not without its disadvantages.  First of all, sentences can have multiple prenex normal forms, with distinct prefixes.  Thus, a sentence could have complexity $\Pi_{1}$ or $\Sigma_{2}$.  This possibility is somewhat troubling for our purposes, but we can resolve the issue by choosing, whenever there are multiple prenex normal forms, the one lowest in the prenex complexity hierarchy.  (If there are multiple minimal prenex normal forms because we choose one arbitrarily.)  By resolving ambiguity this way, we express a preference for simplicity, however defined~\cite{pambuccian1988simplicity}.  Dealing with quantifier alternations becomes psycholigically taxing rather quickly, even with advanced training in logic and mathematics; if there is therefore a prenex normal form that has fewer quantifier alternations than another prenex normal form, is is thus not unreasonable to opt for the simpler one on the theory that, since a human formalizer wrote down the formula, he probably did not have in mind something more complicated.}

{Different prenex normal forms of sentences can influence proof search done with these sentences, but such differences are of no concern here because, upon computing a minimal prenex normal form, no further computation is done with it.  For a discussion, see~\cite{baaz2001normal}.}

{Another disadvantage of prenex normal form is that it can render straightforwardly compehensible sentences into hard-to-understand logical gibberish.  Although this is undoubtedly correct, we needn't be bothered by such rewritings because we don't actually care about the comprehensibility of prenex normal forms, but just their prefixes.}

{Although the language of {\mizar} permits one to write statements using the standard battery of first-order connectives ($\wedge$, $\vee$, $\rightarrow$, $\leftrightarrow$) and quantifiers ($\forall$, $\exists$)~\cite{mizar-in-a-nutshell}, {\mizar} works internally only with $\battery$ (see~\cite{wiedijk-checker} for more details).  The surface forms are rewritten into the more restricted internal form by making use of standard results on interdefinability of connectives and quantifiers for classical first-order logic~\cite{wiedijk-checker}.}

{To more accurately assess the complexity of {\mizar} sentences, we work with the internal (semantic) representations rather than with the surface level layer of sentences that formalizers actually use.  Working with the internal representations is justified because these formulas are what {\mizar} is actually computing with.  Since existential quantifers do not appear in the internal representation, we work with them indirectly as negated universals.  The prenex normal forms of a surface formula $\phi$ are therefore generally different from the prenex normal forms of $\phi$'s internal representation.}

{Due to the austere battery $\battery$ in which all {\mizar} sentences are internally formulated, some theorems have a prenex complexity somewhat greater than they would be if a different battery of connectives were chosen.  For example, in the {\mizar} language one can write formulas such as}

\begin{equation}\label{star}
\forall x (\phi \rightarrow \psi),
\end{equation}

{is represented internally (ignoring differences between $\phi$ and $\psi$ their internal representations) as}

\begin{equation}\label{star-star}
\forall x \neg(\phi \wedge \neg\psi)
\end{equation}

{If $\phi$ is itself a universal formula, say, $\forall y \chi$, then the minimal prenex complexity of (\ref{star}) is no greater than the minimal prenex complexity of}

\[
\forall x \forall y (\chi \rightarrow \psi).
\]

{If one postulates further that $\chi$ and $\psi$ are atomic formulas, then the complexity of (\ref{star}) is $\forall$.  Turning our attention to (\ref{star-star}), we find that the complexity is $\forall\exists$.}

{Going the other way around, prenex normal form can be ``abused'' into being misleadingly low:}

\begin{fact}
{For every theorem $\phi$ of $\ZFC$, there exists an extension by definitions $\ZFC^{\prime}$ of $\ZFC$ and an $0$-ary atomic formula $p$ such that $\proves{\ZFC^{\prime}}{p}$.}
\end{fact}

{There is nothing mysterious going on here: Simply introduce the new $0$-ary predicate symbol $p$ and postulate the definition}
\[
\equivalence{p}{\phi}.
\]

{This amounts essentially to naming a theorem.}

{Although it may be true that we introduce definitions whenever our statements become too complicated, somehow conventional mathematical practice tolerates, or even encourages a mild amount of complexity (quantifier alternation).  Indeed, one could even claim (contra the possibility of trivialization of any theorem via ``atomization''), that some mathematical statements inherently involve some amount of quantifier alternation, such as the notion of continuity of a real function.}


{The answer is that we are interested in taking {\mizar} theorems as-is.  The fact that a human formalizer, with knowledge of logic and mathematics, formalized a theorem in a certain way, suggests that the as-formalized sentence has a mathematical naturalness to it.  That is, after all, one of the central goals of the {\mizar} project.}

\section{Results}

{We computed the minimal prenex normal form complexity for all theorems across the {\mml}\footnote{(We worked with {\mizar} system version 8.1.02 and {\mml} version 5.20.1189.)}.  For us, a ``theorem'' of the {\mml} includes rather more than just those items that are explicitly called \verb+theorem+ in the official {\mizar} syntax.  In addition to official theorems, by ``theorem'' we mean:}

\begin{itemize}
\item {Toplevel propositions (``lemmas''), of various kinds:}
\begin{itemize}
\item {Straightforward statements (usually these are ``article-local'' lemmas that a formalizer prefers not to ``export'', that is, to make it available to any other formalizer}
\item {So-called diffuse reasoning blocks}
\item {Statements introduced by (in effect) existential introduction, (``Consider an $x$ such that\dots'')}
\item {Type-changing statements (statements to the effect that an object ought to be seen as belonging to a type other than what it directly has~\cite{Wiedijk:2007:MST:1792233.1792261}, similar to ``casting'' in programming language terminology)}
\end{itemize}
\item {Properties~\cite{Naumowicz:2004:IMT}, such as the property of $\leq$ being a connected relation on the set of real numbers or the property that inversion on complex numbers is an involution}
\item {Rewrite rules~\cite{Kornilowicz:2013:RRM}}
\item {Identifications}
\item {Compatibility of redefinitions.  (An already-defined notion can be defined in a new way, in which case {\mizar} requires one to prove that the new definiens is equivalent to the old definiens.  One can also define a notion to have a relatively general type, and then redefine it to have a new type, which also requires some justification.)}
\item {Scheme instances  (Schemes themselves, such as the axiom scheme of replacement or the scheme of mathematical induction, are essentially parameterized sets of theorems, and are thus not, strictly speaking, theorems.)}
\item {Existence and uniqueness conditions for new function symbols (these ensure conservativity when extending the language)}
\item {Existence conditions for types (another conservativity condition: in the {\mizar} logic, one can introduce new types only on the condition that they are provably non-empty)}
\item {Definitional theorems.}
\end{itemize}

{Tables~\ref{tab:theorem-data-table} presents our findings, restricted to theorems; Table~\ref{tab:definition-data-table} is restricted just to definitional theorems.}

\begin{table}
\begin{tabular}{l|l|l|}
$n$ & Number of $\Pi_{n}$ sentences & Number of $\Sigma_{n}$ sentences\\
\hline
$0$ & $3254$ & $3254$\\
$1$ & $67252$ & $1374$\\
$2$ & $12596$ & $146$\\
$3$ & $5456$ & $9$\\
$4$ & $567$ & $10$\\
$5$ & $355$ & $1$\\
$6$ & $12$ & $4$\\
$7$ & $14$ & $0$\\
\hline
Total & $89506$ & $4798$
\end{tabular}
\caption{\label{tab:theorem-data-table}Numbers of theorems in the {\mml} by their prenex complexities.}
\end{table}

\begin{table}
\begin{tabular}{l|l|l|}
$n$ & Number of $\Pi_{n}$ sentences & Number of $\Sigma_{n}$ sentences\\
\hline
$0$ & $270$ & $270$\\
$1$ & $4886$ & $0$\\
$2$ & $3541$ & $0$\\
$3$ & $1854$ & $0$\\
$4$ & $71$ & $0$\\
$5$ & $34$ & $0$\\
$6$ & $1$ & $0$\\
\hline
Total & $10657$ & $270$
\end{tabular}
\caption{\label{tab:definition-data-table}Numbers of definitional theorems in the {\mml} by their prenex complexities.}
\end{table}

{$\Pi_{m}$ and $\Sigma_{m}$, restricted to theorems, are empty for all $m > 7$.  $\Pi_{m}$ and $\Sigma_{m}$, restricted to definitions, are empty for all $m > 6$  The counts for $\Pi_{0}$ and $\Sigma_{0}$ are identical in both tables because they are the same set of sentences (that is, they are ``singular propositions'', that is, they have no quantifiers.)}


\section{Discussion}

{Table~\ref{tab:theorem-data-table} shows clearly that the lion's share of {\mizar} theorems belongs to a $\Pi$ hierarchy rather than the $\Sigma$ hierarchy.  This likely accords with standard mathematical practice; it does seem that genuinely existential theorems are uncommon.  (Existential statements that involve parameters are, according to our analysis, not purely existential but rather belong to some level of the $\Pi$ hierarchy because the parameters are implicitly universally quantified at the outermost level of the statement.)}

{In the table of definitions (Table~\ref{tab:definition-data-table}), notice that there are no $\Sigma_{k}$ definitions for $k > 0$.  This is as expected, because when one follows the standard format for extension by definitions~\cite[chapter 8]{suppes1999introduction} a definitional theorem can never be a properly existential formula, that is, a $\Sigma_{k}$ for some $k > 0$.  In other words, the only $k$ for which a definitional theorem belongs to $\Sigma_{k}$ is $k = 0$.  Such definitions are equations that don't involve quantifiers, not even implicitly present ones.  For example,}
\[
0 = \emptyset
\]
{(definition of the natural number $0$ as the empty set) or}
\[
{\arcsin} = ({\sin} \mid \left [ -\pi/2, \pi/2 \right ])^{\smile}
\]
{(definition of $\arcsin$, as a relation, as the converse ($^{\smile}$) of $\sin$ (considered as a relation) composed with the closed interval $[-\pi/2,\pi/2]$), or}
\[
e = \sum_{k \geq 0} \frac{1}{k!}
\]
{(the definition of Euler's constant $e$).}


{The most complex theorems can be divided into a few classes:}

\begin{itemize}
\item {Theorems characterizing of a class of mathematical objects as satisfying a certain list of properties, each of which is a mildly complicated formula.  When such a theorem is put into prenex normal form, it understandably can becomes rather complex (even when we take the least complicated prenex normal form);}
\item {Reprepresentation-style theorems saying that every object in a class has, or is isomorphic to, a certain standard form;}
\item {Existence or uniqueness conditions for functions defined with complicated definientia}
\item {``Inherently'' complicated theorems.}
\end{itemize}

{As an example of the second type, consider the following characterization theorem for finite commutative groups in terms of direct products~\cite[Theorem 34]{GROUP_17.ABS}:}

\begin{quote}
  {Let $G$ be a finite commutative group.  Suppose $|G| > 1$. Then there exists:}
  \begin{itemize}
    \item {a non empty finite set $I$}
    \item {an associative group-like commutative multiplicative magma family $F$ of $I$, and}
    \item {a homomorphism $H_{0}$ from $\Pi F$ to $G$}
  \end{itemize}
  {such that}
  \begin{itemize}
  \item {$I = PrimeFactorization(G)$,}
  \item {for every element $p$ of $I$, $F(p)$ is a subgroup of $G$,}
  \item {$F(p) = (PrimeFactorization(G))(p)$, and}
  \item {for every element $p$ and $q$ of $I$ such that $p \neq q$ we have}
    \begin{itemize}
    \item {$F(p) \cap F(q) = \mathbf{1}_{G}$,}
    \item {$H_{0}$ is bijective, and}
    \item {for every $G$-valued total $I$-defined function $x$ satisfying}
      \begin{quote}
        {for every element $p$ of $I$, $x(p) \in F(p)$}
      \end{quote}
          {we have $x \in \Pi F$ and $H_{0}(x) = \Pi x$.}
    \end{itemize}
  \end{itemize}
\end{quote}

{Without looking up relevant definitions, it is evident that this theorems is quite complex (its complexity is $\Pi_{6}$).}

{As an example of the third kind, consider the definition of the set of free variables of a formula in the language of set theory~\cite{ZF_MODEL.ABS} in which one explicitly enumerates the various kinds of formulas---equations, membership statements (``$x in y$''), negations, conjunctions, quantified statements.  The definition is rather long and explicit; it is perhaps no surprise that the existence condition for function defined by a long definiens, when put into prenex normal form, is rather complicated.}

{The fourth category is undefined; it remains unclear what it means to say that a theorem is inherently complex.  It seems that some theorems in the {\mml} could be simplified (in terms of quantifier complexity) by postulating new definitions.}

{As expected, the complexity of theorems and definitions drops off very quickly.  Most theorems are straightforwardly universal ($\Pi_{1}$), though also as expected there are a fair number of universal-existential ($\Pi_{2}$) statements and even universal-existential-universal ($\Pi_{3}$).  These counts can perhaps be seen as a confirmation of the abstract style of mathematics that developed mainly in the 19th century, and which seems to be especially dominant when formalizing mathematics in set theory.}

\section{Conclusion and future work}

{By viewing the {\mizar} Mathematical Library as a large collection of mathematical theorems formalized in classical first-order set theory, one can adopt an experimental approach to foundations of mathematics and concretely investigate problems.  Such an approach compliments the research methods traditionally used in foundations of mathematics.}

{Our starting point was the question of the sentence complexity of mathematical statements.  The question is, at first blush, somewhat vague.  It becomes more precise when we specify a formal theory which we take as a proxy for ``mathematics''.  In foundations of mathematics, set theory is often (rightly or wrongly) regarded as such a proxy.  Even more precisely, one can single out a specific set theory (e.g., $\ZFC$) and, as far as possible, work only within it.  In this paper, the focus was on $\TG$, in the guise of {\mizar} and its associated library of formalized mathematics, the {\mml}.  We looked at the definitions and theorems across the {\mml} and used prenex normal form as a measure of the complexity of these items.}

{Not unexpectedly, although we gained some insight into our original problem, we learn at the same time something about our methods themselves.  Thus, despite its flaws, we take prenex normal form is a measure of sentence complexity.  The simple notion of quantifier alternation becomes the standard by which to compare formal sentences.  Using this metric, we measued the complexity of the contents of the {\mml} (ignoring proofs) and found that, in accordance with one's informal intuition about mathematics, most mathematical statements involve, at most, relatively few quantifier alternations.  We found a handful of extreme cases.  Do these cases show that there are highly complex corners of mathematics, as judged simply by the complexity of the statements (again, ignoring proofs)?  By investigating extreme cases, we found that high quantifier alternation can usually be accounted for by appealing to the original statement without going through the prenex normal form computation.  Thus, we find that prenex normal form can perhaps overestimate (and thereby mislead us) the complexity of statements by assigning a nearly incomprehensible prefix to certain kinds of relatively benign mathematical statements.}

{The {\mizar} language includes Fraenkel terms such as}

\[
\{ f(x) : \phi(x) \}
\]

{and they occur often among the theorems in the {\mml}.  As with definitions introduced whenever complexity threatens to become too high, Fraenkel terms have the effect of bringing down complexity.  In our experiments, Fraenkel terms were treated like any other term, so they were effectively ignored.  One could redo the experiement by eliminating Fraenkel terms.  Presumably, the complexity of the {\mml} would increase because intuitively one uses Fraenkel terms when one wants to refer to somewhat complicated-to-describe sets in a compact way.}

{One can imagine even further experiments in which not just Fraenkel terms, but in which all definitions are expanded.  This would certaily increase complexity (and one would have to address the ambiguity and nondeterminism inherent in ``expanding'' defined notions), but such an experiment would make progress toward the vision of precisely measuring the ``rock bottom'' complexity of real mathematics formalized in set theory.}

\begin{acknowledgements}
\par{To Andrzej Trybulec.}
\end{acknowledgements}

\bibliographystyle{spmpsci}
\bibliography{alama}

\begin{thebibliography}{10}
\providecommand{\url}[1]{{#1}}
\providecommand{\urlprefix}{URL }
\expandafter\ifx\csname urlstyle\endcsname\relax
  \providecommand{\doi}[1]{DOI~\discretionary{}{}{}#1}\else
  \providecommand{\doi}{DOI~\discretionary{}{}{}\begingroup
  \urlstyle{rm}\Url}\fi

\bibitem{DBLP:journals/jar/Alama13}
Alama, J.: Eliciting implicit assumptions of mizar proofs by property omission.
\newblock J. Autom. Reasoning \textbf{50}(2), 123--133 (2013)

\bibitem{alama-mamane-urban}
Alama, J., Mamane, L., Urban, J.: Dependencies in formal mathematics:
  Applications and extraction for \coq{} and \mizar{}.
\newblock In: J.~Jeuring, J.~Campbell, J.~Carette, G.~Dos~Reis, P.~Sojka,
  M.~Wenzel, V.~Sorge (eds.) Intelligent Computer Mathematics, \emph{Lecture
  Notes in Computer Science}, vol. 7362, pp. 1--16. Springer Berlin /
  Heidelberg (2012).
\newblock \urlprefix\url{http://dx.doi.org/10.1007/978-3-642-31374-5_1}

\bibitem{baaz2001normal}
Baaz, M., Egly, U., Leitsch, A.: Normal form transformations.
\newblock In: A.~Robinson, A.~Voronkov (eds.) Handbook of Automated Reasoning,
  vol.~1, chap.~5, pp. 273--331. Elsevier (2001)

\bibitem{ZF_MODEL.ABS}
Bancerek, G.: Models and satisfiability.
\newblock Formalized Mathematics \textbf{1}({\bf 1}), 191--199 (1990)

\bibitem{davis2001history}
Davis, M.: The early history of automated deduction.
\newblock In: A.~Robinson, A.~Voronkov (eds.) Handbook of Automated Reasoning,
  vol.~1, chap.~1, pp. 3--15. North-Holland (2001)

\bibitem{enderton2001mathematical}
Enderton, H.: A Mathematical Introduction to Logic, 2nd edn.
\newblock Academic Press (2001)

\bibitem{friedman1998complexity}
Friedman, H.: complexity of sentences.
\newblock
  \urlprefix\url{http://www.cs.nyu.edu/pipermail/fom/1998-April/001846.html}

\bibitem{mizar-in-a-nutshell}
Grabowski, A., Kornilowicz, A., Naumowicz, A.: {\mizar} in a nutshell.
\newblock Journal of Formalized Reasoning \textbf{3}(2), 153--245 (2010).
\newblock \urlprefix\url{http://jfr.cib.unibo.it/article/view/1980/1356}

\bibitem{GrabowskiTAC:2012}
Grabowski, A., Schwarzweller, C.: Towards automatically categorizing
  mathematical knowledge.
\newblock In: M.~Ganzha, L.A. Maciaszek, M.~Paprzycki (eds.) Federated
  Conference on Computer Science and Information Systems -- FedCSIS 2012,
  Wroclaw, Poland, 9--12 September 2012, Proceedings, pp. 63--68 (2012)

\bibitem{jech2008axiom}
Jech, T.: The Axiom of Choice.
\newblock Dover (2008).
\newblock Unabridged republication of the original 1973 edition by
  North-Holland.

\bibitem{DBLP:conf/aisc/2012}
Jeuring, J., Campbell, J.A., Carette, J., Reis, G.D., Sojka, P., Wenzel, M.,
  Sorge, V. (eds.): Intelligent Computer Mathematics - 11th International
  Conference, AISC 2012, 19th Symposium, Calculemus 2012, 5th International
  Workshop, DML 2012, 11th International Conference, MKM 2012, Systems and
  Projects, Held as Part of CICM 2012, Bremen, Germany, July 8-13, 2012.
  Proceedings, \emph{Lecture Notes in Computer Science}, vol. 7362. Springer
  (2012)

\bibitem{Kornilowicz:2013:RRM}
Korni{\l}owicz, A.: On rewriting rules in {M}izar.
\newblock Journal of Automated Reasoning \textbf{50}(2), 203--210 (2013).
\newblock \doi{10.1007/s10817-012-9261-6}.
\newblock \urlprefix\url{http://dx.doi.org/10.1007/s10817-012-9261-6}

\bibitem{kunen2011set}
Kunen, K.: Set Theory, \emph{Studies in Logic: Mathematical Logic and
  Foundations}, vol.~34.
\newblock College Publications (2011)

\bibitem{mizar-the-first-30-years}
Matuszewski, R., Rudnicki, P.: Mizar: The first 30 years.
\newblock Mechanized Mathematics and Its Applications \textbf{4}, 3--24 (2005).
\newblock \urlprefix\url{http://mizar.org/people/romat/MatRud2005.pdf}

\bibitem{Naumowicz:2004:IMT}
Naumowicz, A., Byli\'nski, C.: Improving {M}izar texts with properties and
  requirements.
\newblock In: A.~Asperti, G.~Bancerek, A.~Trybulec (eds.) Mathematical
  Knowledge Management, Third International Conference, MKM 2004 Proceedings,
  \emph{MKM'04, Lecture Notes in Computer Science}, vol. 3119, pp. 290--301
  (2004).
\newblock \doi{10.1007/978-3-540-27818-4\_21}.
\newblock \urlprefix\url{http://dx.doi.org/10.1007/978-3-540-27818-4_21}

\bibitem{GROUP_17.ABS}
Okazaki, H., Yamazaki, H., Shidama, Y.: {I}somorphisms of direct products of
  finite commutative groups.
\newblock Formalized Mathematics \textbf{21}({\bf 1}), 65--74 (2013).
\newblock \doi{10.2478/forma-2013-0007}

\bibitem{pambuccian1988simplicity}
Pambuccian, V.: Simplicity.
\newblock Notre Dame Journal of Formal Logic \textbf{29}(3), 396--411 (1988)

\bibitem{Pak:2013:MLE}
P\c{a}k, K.: Methods of lemma extraction in natural deduction proofs.
\newblock Journal of Automated Reasoning \textbf{50}(2), 217--228 (2013).
\newblock \doi{10.1007/s10817-012-9267-0}.
\newblock \urlprefix\url{http://dx.doi.org/10.1007/s10817-012-9267-0}

\bibitem{suppes1999introduction}
Suppes, P.: Introduction to Logic.
\newblock Dover (1999).
\newblock Unabridged republication of the original 1957 edition by van
  Nostrand.

\bibitem{tarski1939wellordered}
Tarski, A.: On well-ordered subsets of any set.
\newblock Fundamenta Mathematicae \textbf{32}, 176--183 (1939)

\bibitem{TARSKI.ABS}
Trybulec, A.: Tarski {G}rothendieck set theory.
\newblock Formalized Mathematics \textbf{1}({\bf 1}), 9--11 (1990)

\bibitem{Trybulec:2013:FMM}
Trybulec, A., Korni{\l}owicz, A., Naumowicz, A., Kuperberg, K.: Formal
  mathematics for mathematicians.
\newblock Journal of Automated Reasoning \textbf{50}(2), 119--121 (2013).
\newblock \doi{10.1007/s10817-012-9268-z}.
\newblock \urlprefix\url{http://dx.doi.org/10.1007/s10817-012-9268-z}

\bibitem{DBLP:conf/mkm/Urban03}
Urban, J.: Translating mizar for first order theorem provers.
\newblock In: A.~Asperti, B.~Buchberger, J.H. Davenport (eds.) MKM,
  \emph{Lecture Notes in Computer Science}, vol. 2594, pp. 203--215. Springer
  (2003)

\bibitem{DBLP:conf/mkm/Urban05}
Urban, J.: {XML}-izing {Mizar}: {Making} semantic processing and presentation
  of {MML} easy.
\newblock In: MKM, pp. 346--360 (2005).
\newblock \urlprefix\url{http://dx.doi.org/10.1007/11618027_23}

\bibitem{DBLP:conf/cade/Urban07}
Urban, J.: Malarea: a metasystem for automated reasoning in large theories.
\newblock In: G.~Sutcliffe, J.~Urban, S.~Schulz (eds.) ESARLT, \emph{CEUR
  Workshop Proceedings}, vol. 257. CEUR-WS.org (2007)

\bibitem{wang1960toward}
Wang, H.: Toward mechanical mathematics.
\newblock IBM Journal pp. 2--22 (1960)

\bibitem{wiedijk-checker}
Wiedijk, F.: {CHECKER}.
\newblock Available online at \url{http://www.cs.ru.nl/~freek/mizar/by.ps.gz}.

\bibitem{Wiedijk:2007:MST:1792233.1792261}
Wiedijk, F.: Mizar's soft type system.
\newblock In: Proceedings of the 20th international conference on Theorem
  proving in higher order logics, TPHOLs'07, pp. 383--399. Springer-Verlag,
  Berlin, Heidelberg (2007).
\newblock \urlprefix\url{http://portal.acm.org/citation.cfm?id=1792233.1792261}

\end{thebibliography}

\end{document}